\documentclass[%
onecolumn
]{mpi2015-cscpreprint}

\usepackage{graphicx}      

\usepackage[american]{babel}

\usepackage{type1cm}        
\usepackage{makeidx}         
\usepackage{graphicx}        
\usepackage{authblk}

\usepackage{newtxtext}       %
\usepackage{newtxmath}       
\usepackage{amsmath,amsfonts}

\usepackage{geometry}
\newcommand{\bA}{{\mathbf A}}
\newcommand{\bB}{{\mathbf B}}
\newcommand{\bC}{{\mathbf C}}
\newcommand{\bD}{{\mathbf D}}

\newcommand{\bF}{{\mathbf F}}
\newcommand{\bG}{{\mathbf G}}
\newcommand{\bY}{{\mathbf Y}}
\newcommand{\bL}{{\mathbf L}}

\newcommand{\bN}{{\mathbf N}}
\newcommand{\bI}{{\mathbf I}}
\newcommand{\bH}{{\mathbf H}}
\newcommand{\bQ}{{\mathbf Q}}
\newcommand{\bW}{{\mathbf W}}
\newcommand{\bT}{{\mathbf T}}
\newcommand{\bZ}{{\mathbf Z}}
\newcommand{\bX}{{\mathbf X}}
\newcommand{\bx}{{\mathbf x}}
\newcommand{\by}{{\mathbf y}}
\newcommand{\bu}{{\mathbf u}}
\newcommand{\bV}{{\mathbf V}}
\newcommand{\bU}{{\mathbf U}}
\newcommand{\bfe}{{\mathbf e}}
\newcommand{\bv}{{\mathbf v}}

\newcommand{\cH}{ {\cal H} }
\newcommand{\Si}{ \boldsymbol{\Sigma} }

\newcommand{\bOmega}{ \boldsymbol{\Omega} }
\newcommand{\bLambda}{\boldsymbol{\Lambda}}
\newcommand{\bGamma}{\boldsymbol{\Gamma}}

\newtheorem{remark}{Remark}[section]

\makeindex             

\begin{document}
	
		\title{Toward fitting structured nonlinear systems by means of dynamic mode decomposition} 
	
	\author[$\ast$]{Ion Victor Gosea}
	\affil[$\ast$]{Max Planck Institute for Dynamics of Complex Techincal Systems, Sandtorstrasse 1, 39106 Magdeburg, Germany, Data-Driven System Reduction and Identification (DRI).\authorcr
		\email{gosea@mpi-magdeburg.mpg.de}, \orcid{0000-0003-3580-4116}}
	
	\author[$\dagger$]{Igor Pontes Duff}
	\affil[$\dagger$]{Max Planck Institute for Dynamics of Complex Techincal Systems, Sandtorstrasse 1, 39106 Magdeburg, Germany, Computational Methods in Systems and Control Theory (CSC).\authorcr
		\email{pontes@mpi-magdeburg.mpg.de}, \orcid{0000-0001-6433-6142}}
	
	\shorttitle{Bilinear $\&$ Quadratic DMD}
	\shortauthor{Ion Victor Gosea, Igor Pontes Duff}
	\shortdate{}
	

%

\abstract{The dynamic mode decomposition (DMD) is a data-driven method used for identifying the dynamics of complex nonlinear systems. It extracts important characteristics of the underlying dynamics using measured time-domain data produced either by means of experiments or by numerical simulations. In the original methodology, the measurements are assumed to be approximately related by a linear operator. Hence, a linear discrete-time system is fitted to the given data. However, often, nonlinear systems modeling physical phenomena have a particular known structure. In this contribution, we propose an identification and reduction method based on the classical DMD approach allowing to fit a structured nonlinear system to the measured data. We mainly focus on two types of nonlinearities: bilinear and quadratic-bilinear. By enforcing this additional structure, more insight into extracting the nonlinear behavior of the original process is gained. Finally, we demonstrate the proposed methodology for different examples, such as the Burgers' equation and the coupled van der Pol oscillators.}

\maketitle


\section{Introduction}

Mathematical models are commonly used to simulate, optimize, and control the behavior of real dynamical processes. A common way to derive those models is to use the first principles, generally leading to a set of ordinary or partial differential equations. For high complex dynamics,  fine discretization leads to high fidelity models, which require numerous equations and variables.  In some situations,  the high model is given as a black box setup, i.e., by solvers that allow the computation of the full model states for a given set of initial conditions and inputs, but does not provide the dynamical system realization. In order to better understand such dynamics processes,  it is often beneficial to construct surrogate models using simulated data.   This justifies the development of identification or data-driven model reduction methods. Indeed,  with the ever-increasing availability of measured/simulated data in different scientific disciplines, the need for incorporating this information in the identification and reduction process has steadily grown. The data-driven model reduction problem consists of determining low-order models from the provided data obtained either by experimentation or numerical simulations.  Methods such as Dynamic Mode Decomposition (DMD) have drawn considerable research endeavors.  

DMD is a data-driven method for analyzing complex systems. The purpose is to learn/ extract the important dynamic characteristics (such as unstable growth modes, resonance, and spectral properties) of the underlying dynamical system by means of measured time-domain data. These can be acquired through experiments in a practical setup or artificially, through numerical simulations (by exciting the system). It was initially proposed in \cite{schmid2010dynamic} in the context of analyzing numerical and experimental fluid mechanics problems. Additionally, it is intrinsically related to the Koopman operator analysis, see \cite{rowley2009spectral, mezic2013analysis}. Since its introduction, several extensions have been proposed in the literature, e.g., the exact DMD \cite{tu2014dynamic}, the extended DMD \cite{chen2012variants} and the higher-order DMD \cite{le2017higher}. Also, in order to address control problems, DMD with control inputs was proposed in \cite{pbk16}, and then extended to the case where outputs are also considered in \cite{annoni2016wind, BHM18}. The reader is refereed to \cite{kutz2016dynamic} for a comprehensive monograph on the topic.

Often, nonlinear systems modeling physical phenomena have a particular known structure, such as bilinear and quadratic terms. In the present work, our primary goal is to embed nonlinear structures in the DMD framework. To this aim,  we propose an identification and data-driven reduction method based on the classical DMD approach allowing to fit a bilinear and quadratic-bilinear structures to the measured data.  The choice to fit such terms is due to the fact most systems with analytical nonlinearities (e.g., rational, trigonometrical, polynomial) can be exactly reformulated as quadratic-bilinear systems \cite{Gu11}. Our work is rooted in the two variants, DMD with control and input-output DMD, and can be considered as an extension of those methodologies.

There exist vast literature on learning nonlinear dynamics from data, and we review the most relevant literature for our work. One approach is the so-called Loewner framework, which enables to construct low-order models from frequency domain data. It was initially proposed in \cite{mayo2007framework}, and later extended to bilinear \cite{AGI16} and quadratic-bilinear case \cite{GA18}. Another approach  is the operator inference, proposed \cite{peherstorfer2016data}.  This approach infers polynomial low-order models as a solution of a least-squares problem based on the initial conditions, inputs, trajectories of the states.  This approach was recently extended to systems with non-polynomials \cite{benner2020operator}. Also, the authors in \cite{qian2019transform} show how the use of  lifting transformations can be beneficial to identify the system. Finally, the approach proposed in \cite{peherstorfer2019sampling}, introduces a method based on operator inference enabling to learn exactly the reduced
models that are traditionally constructed with model reduction. It is worth mentioning that the operator inference   approach \cite{peherstorfer2016data}  can be seen as an extension to DMD for nonlinear systems. Indeed, in the this framework,  the reduced order model are allowed to have polynomial terms on the state and its matrices are obtained by solving a least square problems.  The main difference is that this optimization problem is set using the reduced trajectories as the data (see the introduction of \cite{peherstorfer2016data}  for more details). 

In this work,  we aim at fitting nonlinear model structures using the DMD setup, i.e., by using the full-model trajectories, which is the main difference from \cite{peherstorfer2016data}. Additionally, besides the quadratic structure on the state, we also consider reduced-order models having bilinear structure on the state and input.

The rest of the paper is organized as follows. In Section \ref{sec:DMD},  recalls some  results on the classical DMD, DMD with control and the input-output DMD.  In Section \ref{sec:Ext}   we present the main contribution of the paper, which is the incorporation of bilinear and quadratic-bilinear terms in the DMD setup. Finally, in Section  \ref{sec:Num},  we demonstrate the proposed methodology for different examples, such as the Burgers' equation and the coupled van der Pol oscillators.


\section{Dynamic mode decomposition}
\label{sec:DMD}

In this section, we briefly recall the classical DMD framework \cite{schmid2010dynamic}. To this aim, we analyze time-invariant systems of ordinary differential equations (ODEs) written compactly in a continuous-time setting as follows
\begin{equation}\label{nl_ode}
\dot{\bx}(t) = f(\bx(t)),
\end{equation}
where $\bx(t) \in \mathbb{R}^n$ is the state vector and $f:  \mathbb{R}^n \rightarrow \mathbb{R}^n$ is the system nonlinearity. 

By means of sampling the variable $\bx$ in (\ref{nl_ode}) at uniform intervals of time, we collect a series of vectors $ \bx(t_k)$ for sampling times $t_0, t_1, \ldots, t_m$. For simplicity, denote $\bx_k := \bx(t_k)$.

DMD aims at analyzing the relationship between pairs of measurements from a dynamical system. The measurements $\bx_k$ and $\bx_{k+1}$, as previously introduced, are assumed to be approximately related by a linear operator for all $k \in \{0,1,\ldots,m-1\}$.
\begin{equation} \label{recurrence}
\bx_{k+1} \approx \bA \bx_k,
\end{equation}
where $\bA \in \mathbb{R}^{n \times n}$. This approximation is assumed to hold for all pairs of measurements. Next, group together the sequence of collected snapshots of the discretized state $\bx(t)$ and use the following notations:
\begin{align}\label{snapX}
\bX &= \left[ \begin{array}{cccc}
\bx_0 & \bx_1 &  \ldots & \bx_{m-1}
\end{array} \right] \in \mathbb{R}^{n \times m} , \ \
\bX_s = \left[ \begin{array}{cccc}
\bx_1 & \bx_2 &  \ldots & \bx_{m}
\end{array} \right] \in \mathbb{R}^{n \times m}.
\end{align}
The DMD method is based on finding a best-fit solution of an operator A so that the following relation is (approximately) satisfied
\begin{equation} \label{eq_DMD}
\bX_s = \bA \bX,
\end{equation}
which represents the block-version of equation (\ref{recurrence}). Moreover, the above relation does not need to hold exactly. Previous work has theoretically justified using this approximating operator on data generated by nonlinear dynamical systems. For more details, see \cite{tlrbk14}. A best-fit solution is explicitly given as follows:
\begin{equation} \label{sol_DMD}
\bA = \bX_s \bX^{\dagger},
\end{equation}
where $\bX^\dagger \in \mathbb{R}^{m \times n}$ is the Moore-Penrose inverse of matrix $\bX \in \mathbb{R}^{n \times m}$. In the above statement, by "best-fit" it is meant the solution that minimizes the least-squares error in the Frobenius norm (see \cite{BHM18}). More precisely, the matrix $\bA$ in (\ref{sol_DMD}) is the solution of the following optimization problem
\begin{equation}\label{opt_prob1}
\underset{\hat{\bA} \in \mathbb{R}^{n \times n}}{\text{arg min}}\big{(} \Vert \bX_s - \hat{\bA} \bX \Vert_{\text F} \big{)}.
\end{equation}

The so-called DMD modes are given by the eigenvectors of matrix $\bA$ in (\ref{sol_DMD}), collected in matrix $\bT$ with $\bA = \bT \bLambda \bT^{-1}$. These spatial modes of system (\ref{nl_ode}) are computed at a single frequency and are connected to the Koopman operator, see \cite{mezic2013analysis}.

In this work, we will mainly focus on the construction of the reduced-order model rather than the evaluation of the DMD modes.


\subsection{Dynamic mode decomposition with control (DMDc)}
\label{sec:DMDc_lin}

Dynamic mode decomposition with control (DMDc) was introduced in \cite{pbk16} and it modifies the basic framework characterizing DMD. The novelty is given by including measurements of a control input $u(t) \in \mathbb{R}$. It is hence assumed that the dynamics of the original system of ODEs includes an input dependence, i.e. 
\begin{equation}\label{nl_ode_ctr}
\dot{\bx}(t) = f(\bx(t), u(t)),
\end{equation}
which represents a directs extension of (\ref{nl_ode}). In (\ref{nl_ode_ctr}), it is assumed that $f:  \mathbb{R}^n \times \mathbb{R} \rightarrow \mathbb{R}^n$. Then, continue as in the classical DMD case without control to collect a discretized solution $\bx$ at particular time instances.

In this setup, a trio of measurements are now assumed to be connected. The goal of DMDc is to analyze the relationship between a future state measurement $\bx_{k+1}$ with the current measurement $\bx_k$ and the current control $u_k$.

The motivation for this method is that, understanding the dynamic characteristics of systems that have both internal dynamics and applied external control is of great use for many applications, such as for controller design and sensor placement.

The DMDc method is used to discover the underlying dynamics subject to a driving control input by quantifying its effect to the time-domain measurements corresponding to the underlying dynamical system.

A pair of linear operators represented by matrices $\bA \in \mathbb{R}^{n \times n}$ and $\bB \in \mathbb{R}^{n}$ provides the following dependence for each trio of measurement data snapshots $(\bx_{k+1},\bx_k,\bu_k)$
\begin{equation}\label{eqLIN}
\bx_{k+1} = \bA \bx_k + \bB u_k, \ \ 0 \leq k \leq m-1. 
\end{equation}
Next, denote the sequence of control input snapshots with
\begin{align}\label{snapU}
\bU &= \left[ \begin{array}{cccc}
\bu_0 & \bu_1 &  \ldots & \bu_{m-1}
\end{array} \right] \in \mathbb{R}^{1 \times m}.
\end{align}
The first step is to augment the matrix $\bX$ with the row vector $\bU$ and similarly group together the $\bA$ and $\bB$ matrices by using the notations:
\begin{equation}\label{notLIN}
\bG = [\bA \  \ \bB] \in \mathbb{R}^{n \times (n+1)}, \ \ \ \bOmega = \left[ \begin{array}{c}
\bX \\ \bU
\end{array}  \right] \in \mathbb{R}^{ (n+1) \times m}.
\end{equation}
The matrix $\bG$ introduced above will be referred to as the system matrix since it incorporates the matrices corresponding to the system to be fitted.

By letting the index k vary in the range $\{0,1,\ldots,m-1\}$, one can compactly rewrite the $m$ equations in the following matrix format:
\begin{equation}\label{eq_DMDc_lin}
\bX_s = \bA \bX + \bB \bU = [\bA \  \ \bB] \left[ \begin{array}{c}
\bX \\ \bU
\end{array}  \right] :=  \bG \bOmega.
\end{equation}
Thus, similar to standard DMD, compute a pseudo-inverse and solve for $\bG$ as
\begin{equation}\label{sol_DMDc_lin}
\bG = \bX_s \bOmega^{\dagger} \ \ \Rightarrow \ \ [\bA \  \ \bB] =  \bX_s \left[ \begin{array}{c}
\bX \\ \bU
\end{array}  \right] ^{\dagger}.
\end{equation}
The matrix $\bG \in \mathbb{R}^{n \times (n+1)}$ in (\ref{sol_DMDc_lin}) is actually the solution of the following optimization problem
\begin{equation}\label{opt_prob2}
\underset{\hat{\bG} \in \mathbb{R}^{n \times (n+1)}}{\text{arg min}}\Big{(} \Big{\Vert} \bX_s - \hat{\bG} \left[ \begin{array}{c}
\bX \\ \bU
\end{array}  \right] \Big{\Vert}_{\text F} \Big{)}.
\end{equation}
To explicitly compute the matrix in (\ref{sol_DMDc_lin}), we first find the singular value decomposition (SVD) of the augmented data matrix $\bOmega$ as follows
\begin{equation}\label{svd_DMDc_lin}
\bOmega = \bV \Si \bW^T \approx  \tilde{\bV} \tilde{\Si} \tilde{\bW}^T,
\end{equation}
where the full scale and reduced-order matrices have the following dimensions 
$$
\begin{cases} \bV \in \mathbb{R}^{(n+1) \times (n+1)}, \ \ \Si \in \mathbb{R}^{(n+1) \times m}, \ \ \bV \in \mathbb{R}^{m \times m}, \\ \tilde{\bV} \in \mathbb{R}^{(n+1) \times p}, \ \ \tilde{\Si} \in \mathbb{R}^{p \times p}, \ \ \tilde{\bV} \in \mathbb{R}^{m \times r}.  \end{cases}
$$
The truncation index is denoted with $p$, where $p \leqslant n$. The pseudo-inverse $\bOmega^{\dagger}$ is computed using the matrices from the SVD in (\ref{svd_DMDc_lin}), i. e., as $\bOmega^{\dagger} \approx \tilde{\bW} \tilde{\Si}^{-1} \tilde{\bV}^T$.

By splitting up the matrix $\bV^T$ as $\tilde{\bV}^T = [\tilde{\bV}_1^T \ \  \tilde{\bV}_2^T]$, recover the system matrices as
\begin{equation}
\overline{\bA} = \bX_s \tilde{\bW} \tilde{\Si}^{-1} \tilde{\bV}_1^T, \ \ \overline{\bB} = \bX_s \tilde{\bW} \tilde{\Si}^{-1} \tilde{\bV}_2^T.
\end{equation}
As mentioned in, there is one additional step. By performing another (short) SVD of the matrix $\bX_s$, write
\begin{equation}\label{svd_Xs}
\bX_s  \approx  \hat{\bV} \hat{\Si} \hat{\bW}^T, 
\end{equation}
where $\hat{\bV} \in \mathbb{R}^{(n+1) \times r}, \ \ \hat{\Si} \in \mathbb{R}^{r \times r}, \ \ \hat{\bV} \in \mathbb{R}^{m \times r}$. Note that the two SVDs will likely have different truncation values. The following reduced-order approximations of
$\bA$ and $\bB$ are hence computed as:
\begin{align}
\tilde{\bA} &= \hat{\bV}^T  \overline{\bA} \hat{\bV} = \hat{\bV}^T   \bX_s \tilde{\bW} \tilde{\Si}^{-1} \tilde{\bV}_1^T \hat{\bV} \in \mathbb{R}^{r \times r}, \  \tilde{\bB} = \hat{\bV}^T  \overline{\bB}  = \hat{\bV}^T \bX_s \tilde{\bW} \tilde{\Si}^{-1} \tilde{\bV}_2^T \in \mathbb{R}^{r}.
\end{align}


\subsection{Input-output dynamic mode decomposition}
\label{sec:ioDMD_lin}

In this section we discuss the technique proposed in \cite{annoni2016wind} known as input-output dynamic mode decomposition (ioDMD). This method constructs an input-output reduced-order model and can be viewed an extension of DMDc for the case with observed outputs.
As stated in the original work \cite{annoni2016wind}, this method represents a combination of POD and system identification techniques. The proposed method discussed here is similar in a sense to the algorithms for subspace
state space system identification (N4SID) introduced in \cite{OM94}	 and can be also applied to large-scale systems.

We consider as given a system of ODEs whose dynamics is described by the same equations as in (\ref{nl_ode_ctr}). Additionally, assume that observations are collected in the variable $y(t) \in \mathbb{R}$, as  function of the state variable $\bx$ and of the control $u$, written as
\begin{equation} \label{output_ioDMD_lin}
y(t) = g(\bx(t),u(t)),
\end{equation}
where  $g:  \mathbb{R}^n \times \mathbb{R} \rightarrow \mathbb{R}$.

As before, the next step is to collect snapshots of both variable $\bx(t)$ and of the output $y(t)$ sampled at some positive time instances $t_0,t_1,\ldots t_{m-1}$. Again, for simplicity of the exposition, denote with $y_k := y(t_k)$. 

We enforce the following dependence for each trio of measurement data snapshots given by $(\by_{k},\bx_k,\bu_k)$
\begin{equation}\label{eqLIN_ioDMD}
\by_{k} = \bC \bx_k + \bD u_k, \ \ 0 \leq k \leq m-1. 
\end{equation}

Afterwards, collect the output values in a row vector as follows
\begin{align}\label{snapY}
\bY &= \left[ \begin{array}{cccc}
\by_0 & \by_1 &  \ldots & \by_{m-1}
\end{array} \right] \in \mathbb{R}^{1 \times m}.
\end{align}
The ioDMD method aims at fitting the given set of snapshot measurements collected in matrices $\bX_s, \bX$ and vectors  $\bU$ and $\bY$ to a linear discrete-time system characterized by the following equations
\begin{align}\label{eq_ioDMD_lin}
\begin{split}
\bX_s &= \bA \bX + \bB \bU, \\
\bY &= \bC \bX + \bD \bU,
\end{split}
\end{align}
where as before  $\bA \in \mathbb{R}^{n \times n}$ and $\bB \in \mathbb{R}^{n}$, and also $\bC^T \in \mathbb{R}^{n}, \ \bD \in \mathbb{R}$. Note that the first equation in (\ref{eq_ioDMD_lin}) exactly corresponds to the driving matrix equation of DMDc presented in (\ref{sol_DMDc_lin}). Moreover, write the system of equations in (\ref{eq_ioDMD_lin}) compactly as
\begin{equation} \label{eq_ioDMD_mat_lin}
\left[ \begin{matrix}
\bX_s \\ \bY
\end{matrix} \right] = \left[ \begin{matrix}
\bA & \bB \\ \bC & \bD
\end{matrix} \right] \left[ \begin{matrix}
\bX \\ \bU
\end{matrix} \right].
\end{equation}
Next, we adapt the definition of the system matrix $\bG$ from (\ref{notLIN}) by incorporating an extra line as follows
\begin{equation}\label{notioDMD_lin}
\bG = \left[ \begin{matrix}
\bA & \bB \\ \bC & \bD
\end{matrix} \right] \in \mathbb{R}^{(n+1) \times (n+1)},
\end{equation}
while $\bOmega = \left[ \begin{array}{c}
\bX \\ \bU
\end{array}  \right] \in \mathbb{R}^{ (n+1) \times m}$ is as before. Introduce a new notation that will become useful also in the next sections. It represents an augmentation of the shifted state matrix $\bX_s$ with the output observations vector $\bY$, i.e., 
\begin{equation}
\bGamma = \left[ \begin{array}{c}
\bX_s \\ \bY
\end{array}  \right] \in \mathbb{R}^{ (n+1) \times m}.
\end{equation}

Again, the solution of equation (\ref{eq_ioDMD_mat_lin}) will be computed as a best-fit type of approach. Hence, similarly to the DMDc case, recover the matrices $\bA, \bB$, $\bC$, and $\bD$ by computing the pseudo-inverse of matrix $\bOmega$ and writing
\begin{equation}\label{sol_ioDMD_lin}
\bG = \bGamma \bOmega^{\dagger} \ \ \Rightarrow \ \ \left[ \begin{matrix}
\bA & \bB \\ \bC & \bD
\end{matrix} \right] =  \left[ \begin{matrix}
\bX_s \\ \bY
\end{matrix} \right] \left[ \begin{array}{c}
\bX \\ \bU
\end{array}  \right] ^{\dagger}.
\end{equation} 
The matrix $\bG \in \mathbb{R}^{(n+1) \times (n+1)}$ in (\ref{sol_DMDc_lin}) is actually the solution of the following optimization problem
\begin{equation}\label{opt_prob3}
\underset{\hat{\bG} \in \mathbb{R}^{(n+1) \times (n+1)}}{\text{arg min}}\Big{(} \Big{\Vert} \left[ \begin{array}{c}
\bX_s \\ \bY
\end{array}  \right] - \hat{\bG} \left[ \begin{array}{c}
\bX \\ \bU
\end{array}  \right] \Big{\Vert}_{\text F} \Big{)}.
\end{equation}
Similarly to the procedure covered in Section \ref{sec:DMDc_lin}, one could further lower the dimension of the recovered system matrices by employing an additional SVD of the matrix $\bGamma$, as was done in (\ref{svd_Xs}).

  
  \section{The proposed extensions}\label{sec:Ext}
  
  In this section, we present the main contribution of the paper. We propose extensions of the methods previously introduced in Sections \ref{sec:DMDc_lin} and \ref{sec:ioDMD_lin}, e.g. DMDc and, respectively ioDMD to fit nonlinear structured systems. More specifically, the discrete-time models that are fitted using these procedures will no longer be linear as in (\ref{eq_ioDMD_lin}); the new models will contain nonlinear (bilinear or quadratic) terms.
  
  \subsection{Bilinear Systems}
  \label{sec:BIL_DMD}
  
  Bilinear systems are a class of mildly nonlinear systems for which the nonlinearity is given by the product between the state variable and the control input. More exactly, the characterizing system of ODEs is written as in (\ref{nl_ode_ctr}) but for a specific choice of mapping $f$, i.e.  $f(\bx,\bu) = \bA \bx+ \bN \bx u + \bB u$. Additionally, assume that the observed output $y$ depends linearly on the state $\bx$. Hence, it what follows, we will make use of the following description of bilinear systems (with a single input and a single output) 
  \begin{align}\label{bil_ode_ctr}
  \begin{split}
  \dot{\bx}(t) &= \bA \bx(t) + \bN \bx(t) u(t) + \bB   u(t), \\
  y(t) &= \bC \bx(t),
  \end{split}
  \end{align}
  where the matrix $\bN \in \mathbb{R}^{n \times n}$ scales the product of the state variable $\bx$ with the control input $u$. 
  In practice, bilinear control systems are used to approximate nonlinear systems with more general, analytic nonlinearities. This procedure is known as Carleman's linearization; for more details see \cite{RU81}.\\
  
  Bilinear systems are a class of nonlinear systems that received considerable attention in the last four or five decades. Contributions that range from  realization theory in \cite{IR73}, classical system identification in \cite{Do90} or to subspace identification in \cite{FMO99}. In more recent years (last two decades), model order reduction of bilinear systems (in both continuous and discrete-time domain) was extensively studied with contributions covering balanced truncation in \cite{ZLHY03}, Krylov subspace methods in \cite{BD10}, interpolation-based $\cH_2$ method in \cite{BBD11, BB12}, or data-driven Loewner approach in \cite{AGI16,morKarGA20}.
  
  \subsubsection{The general procedure}
  \label{sec:BIL_DMD_GEN}
  
  We start by collecting snapshots of the state $\bx$ for multiple time instances $t_k$. We enforce that the snapshot $\bx_{k+1}$ at time $t_{k+1}$ depends on the snapshot $\bx_k$ in the following way
  \begin{equation}\label{eqBILIN}
  \bx_{k+1} = \bA \bx_{k} + \bN \bx_{k} u_{k} + \bB u_{k}, \ \ \text{for} \ \  0 \leq k \leq m-1. 
  \end{equation}
  We denote the sequence of state and input snapshots as in (\ref{snapX}) and in (\ref{snapU}). Again, by varying the index $k$ in the interval $\{1,2,\ldots,m-1\}$, one can compactly rewrite the $m-1$ equations in the following matrix format:
  \begin{equation}\label{Xs_DMDc_bilin}
  \bX_s = \bA \bX + \bN \bX \bU_{\bD}+ \bB \bU,
  \end{equation}
  where $\bU_\bD = \text{diag}(u_0,u_1,\ldots,u_{m-1}) \in \mathbb{R}^{m \times m}$. One can hence write $\bU = \bL \bU_\bD $, with  $\bL = [1 \ 1 \ \ldots \ 1] \in \mathbb{R}^{1 \times m}$ and then introduce the matrix $\bZ \in \mathbb{R}^{(n+1) \times m}$ as
  \begin{equation}\label{notZ_DMDc_bil}
  \bZ = \left[ \begin{array}{c} \bL \bU_\bD \\ \bX \bU_\bD
  \end{array} \right] =  \left[ \begin{array}{c} \bL  \\ \bX 
  \end{array} \right] \bU_\bD.
  \end{equation}
  The next step is to augment the matrix $\bX$ with matrix $\bZ$ and denote this new matrix with $\bOmega  \in \mathbb{R}^{ (2n+1) \times m}$ as an extension of the matrix previously introduced in (\ref{notLIN}), i.e.,
  \begin{equation}\label{notBILIN}
  \bOmega = \left[ \begin{array}{c}
  \bX \\ \bZ
  \end{array}  \right].
  \end{equation}
  For the case in which we extend the DMDc method in Section \ref{sec:DMDc_lin} to fitting bilinear dynamics (no output observations), we propose a slightly different definition for the matrix $\bG$. We hence append the matrix $\bN$ to the originally introduced system matrix in (\ref{notLIN}). Then, 
  equation (\ref{Xs_DMDc_bilin}) can be written in a factorized way as $\bGamma = \bG \bOmega$,
  where the matrices for this particular setup are as follows
  \begin{equation}\label{notGBILIN}
  \bG = \left[ \begin{matrix} \bA &  \bB & \bN \end{matrix} \right] \in \mathbb{R}^{n \times (2n+1)}, \ \ \bGamma = \bX_s.
  \end{equation}
  Alternatively, for the case where output observations $y_k$ are also available, we enforce a special bilinear dependence for each trio of measurement data snapshots as
  \begin{equation}\label{eqBILIN_ioDMD}
  \by_{k} = \bC \bx_k + \bF \bx_k u_k + \bD u_k, \ \ 0 \leq k \leq m-1,
  \end{equation}
  where $\bF^T \in \mathbb{R}^{n}$. Note that (\ref{eqBILIN_ioDMD}) represents a natural extension of the relation imposed in (\ref{eqLIN_ioDMD}). There, fitting a linear structure is instead enforced. 
  
  Afterwards, we collect the equations in (\ref{eqBILIN_ioDMD}) for each index $k$ and hence write
  \begin{equation}\label{Y_DMDc_bilin}
  \bY = \bC \bX + \bF \bX \bU_{\bD}+ \bD \bU,
  \end{equation}
  with the same notations as in (\ref{notZ_DMDc_bil}). Then, by combining (\ref{Xs_DMDc_bilin}) and (\ref{Y_DMDc_bilin}), we can write all snapshot matrix quantities in the following structured equalities
  \begin{align}
  \begin{split}
  \bX_s = \bA \bX + \bN \bX \bU_{\bD}+ \bB \bU, \\
  \bY = \bC \bX + \bF \bX \bU_{\bD}+ \bD \bU.
  \end{split}
  \end{align}
  This system of equations can be then written in a factorized way as before, i.e., $\bGamma = \bG \bOmega$, where the matrices for this particular setup are given below
  \begin{equation}\label{not_G_GAM_BILIN}
  \bG = \left[ \begin{matrix} \bA &  \bB & \bN \\ \bC & \bD & \bF \end{matrix} \right] \in \mathbb{R}^{(n+1) \times (2n+1)}, \ \ \bGamma = \left[ \begin{matrix}
  \bX_s \\ \bY
  \end{matrix} \right].
  \end{equation}
  
  Finally, the last step is to recover the matrix $\bG$ and split it block-wise in order to put together a system realization. Consequently, this all boils down to solving the equation $\bGamma = \bG \bOmega$ (in either of the two cases, with or without output observations included). More precisely, the objective matrix $\bG \in \mathbb{R}^{(n+1) \times (2n+1)}$ in (\ref{not_G_GAM_BILIN}) is the solution of the following optimization problem
  \begin{equation}\label{opt_prob4}
  \underset{\hat{\bG} \in \mathbb{R}^{(n+1) \times (2n+1)}}{\text{arg min}}\Big{(} \Big{\Vert} \left[ \begin{array}{c}
  \bX_s \\ \bY
  \end{array}  \right] - \hat{\bG} \left[ \begin{array}{c}
  \bX \\ \bZ
  \end{array}  \right] \Big{\Vert}_{\text F} \Big{)} \Leftrightarrow  \underset{\hat{\bG} \in \mathbb{R}^{(n+1) \times (2n+1)}}{\text{arg min}}\Big{(} \big{\Vert} \bGamma - \hat{\bG} \bOmega \big{\Vert}_{\text F} \Big{)}
  \end{equation}
  As shown in the previous sections, solving for $\bG$ in (\ref{opt_prob4}) involves computing the pseudo-inverse of matrix $\bOmega\in \mathbb{R}^{ (2n+1) \times m}$ from (\ref{notBILIN}). More precisely, we write the solution as
  \begin{equation}
  \bG = \bGamma \bOmega^{\dagger}.
  \end{equation}
  
  \begin{remark}
  	Note that the observation map $g$ corresponding to the original dynamical system, as introduced 
  	in (\ref{output_ioDMD_lin}), need not have a bilinear structure as in (\ref{eqBILIN_ioDMD}). It could include more complex nonlinearities or could even be linear. In the later case, the recovered matrix $\bF$ will typically have a low norm.
  \end{remark}
  
  \subsubsection{Computation of the reduced-order matrices}
  
  In this section we present specific/practical details for retrieving the system matrices in the case of the proposed procedure in Section \ref{sec:BIL_DMD_GEN}. We solve the equation $\bGamma = \bG \bOmega$ for which the matrices are given as in (\ref{not_G_GAM_BILIN}), i.e. the case containing output observations. We compute an SVD of the augmented data matrix $\bOmega$ giving
  \begin{equation}
  \bOmega = \bV \Si \bW^T \approx  \tilde{\bV} \tilde{\Si} \tilde{\bW}^T,
  \end{equation}
  where the full-scale and reduced-scale matrices derived from SVD are as follows
  $$
  \begin{cases} \bV \in \mathbb{R}^{(2n+1) \times (2n+1)}, \ \ \Si \in \mathbb{R}^{(2n+1) \times (m-1)}, \ \ \bV \in \mathbb{R}^{(m-1) \times (m-1)}, \\ \tilde{\bV} \in \mathbb{R}^{(2n+1) \times p}, \ \ \tilde{\Si} \in \mathbb{R}^{p \times p}, \ \ \tilde{\bV} \in \mathbb{R}^{(m-1) \times p}.  \end{cases}
  $$
  The truncation index is denoted with $p$, where $p \leqslant n$.The computation of the pseudo-inverse $\bOmega^{\dagger}$ is done by the SVD approach, i.e., $
  \bOmega^{\dagger} \approx \tilde{\bW} \tilde{\Si}^{-1} \tilde{\bV}^T.$ By splitting up the $\bV^T$ matrix as $\tilde{\bV}^T = [\tilde{\bV}_1^T \ \  \tilde{\bV}_2^T \ \ \tilde{\bV}_3^T]$, one can recover the system matrices as 
  \begin{align} \label{sys_mat_bil}
  \begin{split}
  \overline{\bA} &= \bX_s \tilde{\bW} \tilde{\Si}^{-1} \tilde{\bV}_1^T, \ \ \overline{\bB} = \bX_s \tilde{\bW} \tilde{\Si}^{-1} \tilde{\bV}_2^T, \ \
  \overline{\bN} = \bX_s \tilde{\bW} \tilde{\Si}^{-1} \tilde{\bV}_3^T, \\
  \overline{\bC} &= \bY \tilde{\bW} \tilde{\Si}^{-1} \tilde{\bV}_1^T, \ \ \overline{\bD} = \bY \tilde{\bW} \tilde{\Si}^{-1} \tilde{\bV}_2^T, \ \
  \overline{\bF} = \bY \tilde{\bW} \tilde{\Si}^{-1} \tilde{\bV}_3^T.
  \end{split}
  \end{align}
  By performing another (short) SVD for the matrix $\bX_s$, we can write
  \begin{equation}
  \bX_s  \approx  \hat{\bV} \tilde{\Si} \hat{\bW}^T, 
  \end{equation}
  where $\hat{\bV} \in \mathbb{R}^{(n+1) \times r}, \ \ \hat{\Si} \in \mathbb{R}^{r \times r}, \ \ \hat{\bW} \in \mathbb{R}^{(m-1) \times p}$. Note that the two SVDs could have different truncation values denoted with p and r. Using the transformation $\bx = \hat{\bV} \tilde{\bx}$, the following reduced-order matrices can be computed:
  \begin{align} \label{sys_mat_red_bilin}
  \begin{split}
  \tilde{\bA} &= \hat{\bV}^T  \overline{\bA} \hat{\bV} \in \mathbb{R}^{r \times r} , \ \ \
  \tilde{\bB}  = \hat{\bV}^T  \overline{\bB}  \in \mathbb{R}^{r}, \ \ \
  \tilde{\bN} = \hat{\bV}^T  \overline{\bN} \hat{\bV}  \in \mathbb{R}^{r \times r}, \\
  \tilde{\bC} &=  \overline{\bC} \hat{\bV} \in \mathbb{R}^{1 \times r} , \ \ \
  \tilde{\bD}  =   \overline{\bD}  \in \mathbb{R}, \ \ \
  \tilde{\bF} = \hat{\bV}^T  \overline{\bF} \hat{\bV}  \in \mathbb{R}^{1 \times r}.
  \end{split}
  \end{align}
  
  \subsubsection{Conversions between discrete-time and continuous-time representations}
  
  The DMD-type approaches available in the literature identify continuous-time systems by means of linear discrete-time models. In this contribution, we make use of the same philosophy, in the sense that the models fitted are discrete-time. We extend the DMDc and ioDMD approaches by allowing bilinear or quadratic terms to appear in the these models as well.
  
  As also mentioned in \cite{BHM18}, one can compute a continuous-time
  model that represents a first-order approximation of the discrete time model obtained by DMD-type approaches.
  
  Assume that we are in the bilinear setting presented in Section \ref{sec:BIL_DMD} and that we already have computed a reduced-order discrete-time model given by matrices $\{\tilde{\bA},\tilde{\bB},\tilde{\bN},\tilde{\bC},\tilde{\bD},\tilde{\bF}\}$, i.e., following the explicit derivations in (\ref{sys_mat_red_bilin}). Then, a continuous-time model 
  $\{\hat{\bA},\hat{\bB},\hat{\bN},\hat{\bC},\hat{\bD},\hat{\bF}\}$ can also be derived. By assuming that the standard first-order Euler method was used for simulating the original system (with a small enough time step-size $0 < \Delta_t \ll 1$), we can write that
  \begin{align}\label{conversion}
  \begin{split}
  \bx_{k+1} &= \bx_k + \Delta_t (\hat{\bA} \bx_k+ \hat{\bB} \bu_k+\hat{\bN} \bx_k \bu_k) \Rightarrow \\ \tilde{\bA} \bx_k+ \tilde{\bB} \bu_k + \tilde{\bN} \bx_k \bu_k &= \bx_k  +\Delta_t \hat{\bA} \bx_k+ \Delta_t \hat{\bB} \bu_k + \Delta_t \hat{\bN} \bx_k \bu_k \Rightarrow \\
  & \begin{cases} \hat{\bA}  =  \Delta_t^{-1}(\tilde{\bA}-\bI), \ \ \  \hat{\bB} =  \Delta_t^{-1} \hat{\bB}, \ \ \ \hat{\bN} =  \Delta_t^{-1} \hat{\bN}, \\ 
  \hat{\bC}  = \tilde{\bC}, \ \ \  \hat{\bD}  = \tilde{\bD}, \ \ \ \hat{\bF}  = \tilde{\bF}. 
  \end{cases}
  \end{split}
  \end{align}
  Observe that for the ioDMD-type of approaches, the feed-through terms that appear in the output-state equation are the same in both discrete and continuous representations.
  
  \subsection{Quadratic-Bilinear Systems}
  \label{sec:QB_DMD}
  
  Next, we extend the method in Section \ref{sec:DMDc_lin} for fitting another class of nonlinear systems, i.e. quadratic-bilinear (QB) systems. Additional to the bilinear terms that enter the differential equations, we assume that quadratic terms are also present. More precisely, the system of the ODEs is written as in (\ref{nl_ode_ctr}) but for a specific choice of nonlinear mapping $f$, i.e.,
  $$
  f(\bx,\bu) = \bA \bx+ \bQ (\bx \otimes \bx) + \bN \bx u + \bB u,
  $$
  where '$\otimes$' denotes the Kronecker product, the matrix $\bQ \in \mathbb{R}^{n \times n^2}$ scales the product of the state $\bx$ with itself, and $\bN \in \mathbb{R}^{n \times n}$ is as in Section \ref{sec:BIL_DMD}.
  
  Quadratic-bilinear systems appear in many applications for which the original system of ODEs inherently has the required quadratic structure. For example, after semi-discretizing the Burgers' or Navier-Stokes equations in the spatial domain, one obtains a system of differential equations with quadratic nonlinearities (and also with bilinear terms). Moreover, many smooth analytic nonlinear systems, that contain combinations of nonlinearities such as, e.g. exponential, trigonometric, polynomial functions, etc. can
  be equivalently rewritten as QB systems. This is performed by employing so-called lifting techniques. More exactly, one needs to introduce new state variables in order to simplify the nonlinearities  and hence derive new differential equations corresponding to these variables. Model order reduction of QB systems was a topic of interest in the last years with contributions ranging from projection-based approaches in \cite{Gu11,BB15} to optimal $\mathcal{H}_2$-based approximation in \cite{BGG18}, or data-driven approaches in the Loewner framework in \cite{GA18,AGH18}.
  
  Similarly to the procedure described in Section \ref{sec:BIL_DMD}, we enforce that the snapshot $\bx_{k+1}$ at time $t_{k+1}$ depends on the snapshot $\bx_k$ in the following way
  \begin{equation}\label{eqQB}
  \bx_{k+1} = \bA \bx_{k} + \bQ (\bx_{k} \otimes \bx_{k}) + \bN \bx_{k} u_{k} + \bB u_{k}, \ \ \text{for} \ \  0 \leq k \leq m-1. 
  \end{equation}
  Next, by varying the $k$ in the range $\{1,2,\ldots,m-1\}$,  compactly rewrite the $m$ equations in (\ref{eqQB}) in the following matrix format:
  \begin{equation}\label{eq1Q1}
  \bX_s = \bA \bX + \bQ (\bX \otimes \bX) \bH + \bN \bX \bU_{\bD}+ \bB \bU,
  \end{equation}
  with $\bU_\bD = \text{diag}(u_0,u_1,\ldots,u_{m-1}) \in \mathbb{R}^{m \times m}$ and $\bH = \left[ \begin{array}{cccc}
  \bfe_1 \otimes \bfe_1 & \bfe_2 \otimes \bfe_2 & \ldots & \bfe_m \otimes \bfe_m
  \end{array} \right] \in \mathbb{R}^{m^2 \times m}.$ Here, $\bfe_k$ is the unit vector of length $n$ that contains the 1 on position $k$. Additionally, we introduce the matrix $\bT$ that depends on the state matrix $\bX$ as
  $$
  \bT = \left[ \begin{array}{cccc}
  \bX_1 \otimes \bX_1 & \bX_2 \otimes \bX_2 & \ldots & \bX_{m} \otimes \bX_{m}
  \end{array} \right] \in \mathbb{R}^{n^2 \times m}.
  $$
  Note that the equality holds as follows $\bT = (\bX \otimes \bX) \bH$. Next, we augment the matrix $\bX$ with both matrices $\bZ$ and  $\bT$ and group together the matrices $\bA, \ \bQ, \ \bN$ and $\bB$ by using the notations:
  \begin{equation}\label{notQB}
  \bG = [\bA \ \bB \  \bN \ \bQ] \in \mathbb{R}^{n \times (n^2+2n+1)}, \ \  \bOmega = \left[ \begin{array}{c}
  \bX \\ \bZ \\ \bT
  \end{array}  \right] \in \mathbb{R}^{ (n^2+2n+1) \times m}, \ \ \bGamma = \bX_s.
  \end{equation}
  Hence, by using the above notations, rewrite the equation (\ref{eq1Q1}) as follows $
  \bGamma = \bG \bOmega.$
  
  More precisely, the objective matrix $\bG \in \mathbb{R}^{n \times (n^2+2n+1)}$ in (\ref{notQB}) is the solution of the following optimization problem
  \begin{equation}\label{opt_prob5}
  \underset{\hat{\bG} \in \mathbb{R}^{n \times (n^2+2n+1)}}{\text{arg min}}\Big{(} \big{\Vert} \bGamma - \hat{\bG} \bOmega \big{\Vert}_{\text F} \Big{)}.
  \end{equation}
  Thus, one can recover the matrix $\bG$ by solving an optimization problem, e.g., the one given in (\ref{opt_prob5}). This is explicitly done by computing the Moore-Pseudo pseudo-inverse of matrix $\bOmega \in \mathbb{R}^{ (n^2+2n+1) \times m}$, and then writing  $\bG = \bGamma \bOmega^{\dagger}$.
  
  As previously shown in Section \ref{sec:BIL_DMD}, we can again adapt the procedure for fitting QB systems in the ioDMD format by involving output observations measurements $y_k$. The procedure for quadratic-bilinear systems is similar to that for bilinear systems and we prefer to skip the exact description to avoid duplication. For more details, see the derivation in Section \ref{Appendix_QB}.
  
  \begin{remark}
  	Note that the Kronecker product of the vector $\bx \in \mathbf{R}^n$ with itself, i.e., $\bx^{(2)}  = \bx \otimes \bx$ has indeed duplicate components. For $n=2$, one can write 
  	$$
  	\bx^{(2)} = \left[ \begin{matrix}
  	x_1^2 \ & x_1 x_2 \ & x_2 x_1 \ & x_2^2
  	\end{matrix} \right]^T.
  	$$
  	Thus, since matrix $\bG$ is explicitly written in terms of $\bQ$ as in (\ref{notQB}), the linear system of equations $\bGamma = \bG \bOmega$ does not have an unique solution. By using the Moore-Penrose inverse, one implicitly regularizes the least-squares problem in (\ref{opt_prob5}). Additionally,  note that using a different least-squares solver (with or without regularization) could indeed produce a different result.  
  \end{remark}
  
  \begin{remark}
  	It is to be noted that the operator inference procedure avoids the non-uniqueness by accounting for duplicates in the vector $\bx \otimes \bx$. This is done by introducing a special Kronecker product for which the duplicate terms are removed. For more details, we refer the reader to Section 2.3 from \cite{benner2020operator}.
  \end{remark}	
  
  
  \section{Numerical Experiments}\label{sec:Num}
  
  \subsection{The viscous Burgers' equation}
  
  Consider the partial differential viscous Burgers' equation:
  \begin{equation*}
  \small
  \frac{\partial v(x,t)}{\partial t}+v(x,t) \frac{\partial v(x,t)}{\partial x} =  \nu 
  \frac{\partial^2 v(x,t)}{\partial x^2},\  \ \ \ (x,t) \in (0,L) \times (0,T)\;, 
  \normalsize
  \end{equation*}
  with i.c. $v(x,0) = 0, \ x \in [0,L],\ v(0,t) = u(t),\  v(L,t) = 0, \ t \geqslant 0$. The viscosity parameter is denoted with $\nu$.
  
  The Burgers' equation has a convective term, an unsteady term and a viscous term; it can be viewed as a simplification of the Navier-Stokes equations.
  
  By means of semi-discretization in the space domain, one can obtain the following nonlinear (quadratic) model (see \cite{BB15}) described by the following system of ODEs
  \begin{equation} \label{QB_Burgers}
  \dot{v}_k = \begin{cases} -\frac{1}{2h} v_1v_2+\frac{\nu}{h^2}(v_2-2v_1) + (\frac{v_1}{2h}+\frac{\nu}{h^2})u, \ \ k=1, \\[1mm] -\frac{v_k}{2h}(v_{k+1}-v_{k-1})+\frac{\nu}{h^2}(v_{k+1}-2v_k+v_{k-1}), \ \ 2 \leqslant k \leqslant n_0-1, \\[1mm] -\frac{1}{2h} v_nv_{n-1}+\frac{\nu}{h^2}(-2v_n +2v_{n-1}), \ \ k=n_0. \end{cases} 
  \end{equation}
  Next, by means of the Carleman linearization procedure in \cite{RU81}, one can approximate the above nonlinear system of order $n_0$ with a bilinear system of order $n = n_0^2+n_0$. The procedure is as follows: let $\bv = \left[ v_1 \  v_2 \ \ldots \ v_n \right]^T$ be the original state variable in (\ref{QB_Burgers}). Then, introduce the augmented state variable $\bx =  \left[ \begin{array}{c} \bv \\ \bv \otimes \bv \end{array} \right] \in \mathbb{R}^{n_0^2+n_0}$ corresponding to the system  described by the following equations
  \begin{align} \label{Bil_Burgers}
  \begin{split}
  \dot{\bx} &= \bA \bx + \bN \bx u+\bB u, \\
  \by &= \bC \bx.
  \end{split}
  \end{align}
  The continuous-time bilinear model in (\ref{Bil_Burgers}) is going to be used in following the numerical experiments.
  
  Start by choosing the viscosity parameter to be $\nu = 0.01$. Then, choose $n_0 = 40$ as the dimension of the original discretization and hence, the bilinear system in (\ref{Bil_Burgers}) is of order $n = 1640$. Perform a time-domain simulation of this system  by approximating the derivative as follows $\dot{\bx}(t_k) \approx \frac{\bx(t_{k+1})- \bx(t_{k})}{t_{k+1}-t_k} = \frac{\bx_{k+1}-\bx_{k}}{\Delta_t}$. We use as time step $\delta_t = 10^{-3}$ and the time horizon to be $[0,10]s$. The control input is chosen to be $u(t) = 0.5 \cos(10t)e^{-0.3 t}$.
  
\begin{figure}[!htb]
  	\vspace*{-2ex}
  	\begin{center}
  		\includegraphics[width = 0.8\textwidth]{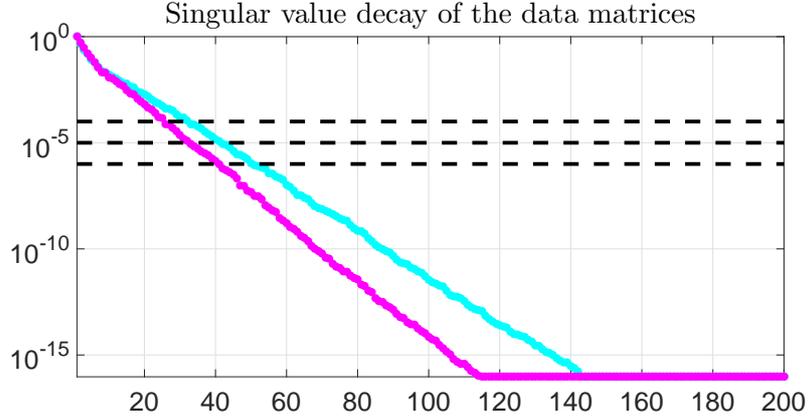} 
  		
  		\caption{The normalized first 200 singular values of matrices $\bOmega$ (with cyan)  and $\bGamma$ (with magenta). The three dotted black lines correspond to the three tolerance levels chosen for $\tau_r$.}
  		\label{fig:GP1}
  	\end{center}
  	\vspace{-3mm}
  \end{figure}

  Hence, collect $10^4$ snapshots of the trio $(\bx_k,u_k,y_k)$ that are arranged in the required matrix format as presented in the previous sections. The first step is to perform an SVD for the matrix $\bOmega \in \mathbb{R}^{3281 \times 10^4}$. The first 200 normalized singular values are presented in Figure\;\ref{fig:GP1}. Choose the tolerance value $\tau_p = 10^{-10}$ which corresponds to a truncation order of $p = 86$ (for computing the pseudo-inverse of matrix $\bOmega$). On the same plot in Figure\;\ref{fig:GP1} we also display the normalized singular values of matrix $\bGamma \in \mathbb{R}^{1641 \times 10^4}$. Note that machine precision is reached at the $112^{\text th}$ singular value. We select three tolerance values $\tau_r \in \{ 10^{-4}, 10^{-5}, 10^{-6} \}$ for truncating matrices obtained from the SVD of $\bGamma$. 
  
  In what follows, we compute reduced-order discrete-time models that have dimension $r$, as in (\ref{sys_mat_red_bilin}). Next, these models are converted using (\ref{conversion}) into a continuous-time model. 
  
  \subsubsection{Experiment 1 - validating the trained models} \label{sec:Burgers_Exp1}
  
  In this first setup, we perform time-domain simulations of the reduced-order models for the same conditions as in the training stage, i.e. in the time horizon $[0,10]s$ and by using the control input $u(t) = 0.5 \cos(10t)e^{-0.3 t}$. Hence, we are validating the trained models on the training data.
  
\begin{figure}[!htb]
  	\vspace*{-2ex}
  	\begin{center}
  		\includegraphics[width = 0.9\textwidth]{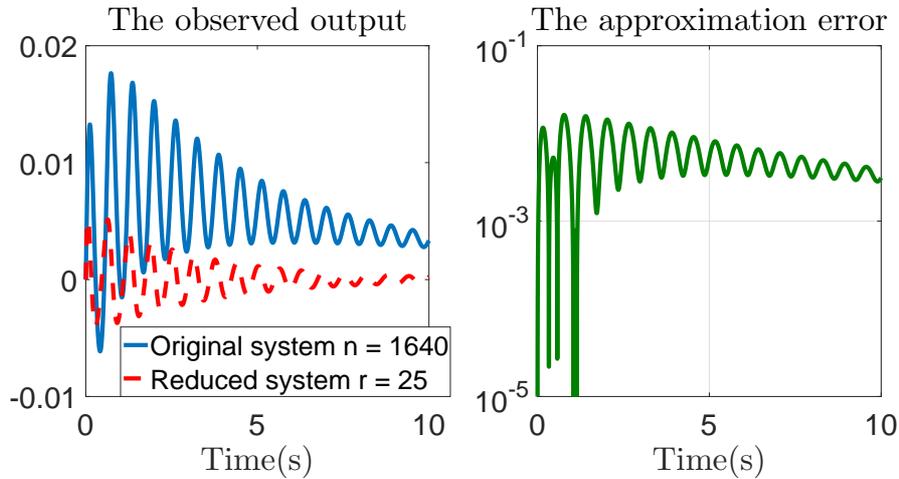} 
  		
  		\caption{Left plot: the observed outputs; 
  			right plot: the corresponding approximation error.}
  		\label{fig:GP2}
  	\end{center}
  	\vspace{-3mm}
  \end{figure}

  Start by choosing the first tolerance value, e.g. $\tau_r = 10^{-4}$. This corresponds to a truncation value of $r=25$. We compute term $\hat{\bD} = 1.1744\text e-14$ and a also bilinear feed-through term with $\Vert \hat{\bF} \Vert_2 = 6.7734\text e-04$. We simulate both the original large-scale bilinear system and the reduced-order system. The results are presented in  Figure\;\ref{fig:GP2}. Note that the observed output curves deviate substantially. One way to improve this behavior is to decrease the tolerance value.
  
\begin{figure}[!htb]
  	\vspace*{-2ex}
  	\begin{center}
  		\includegraphics[width = 0.9 \textwidth]{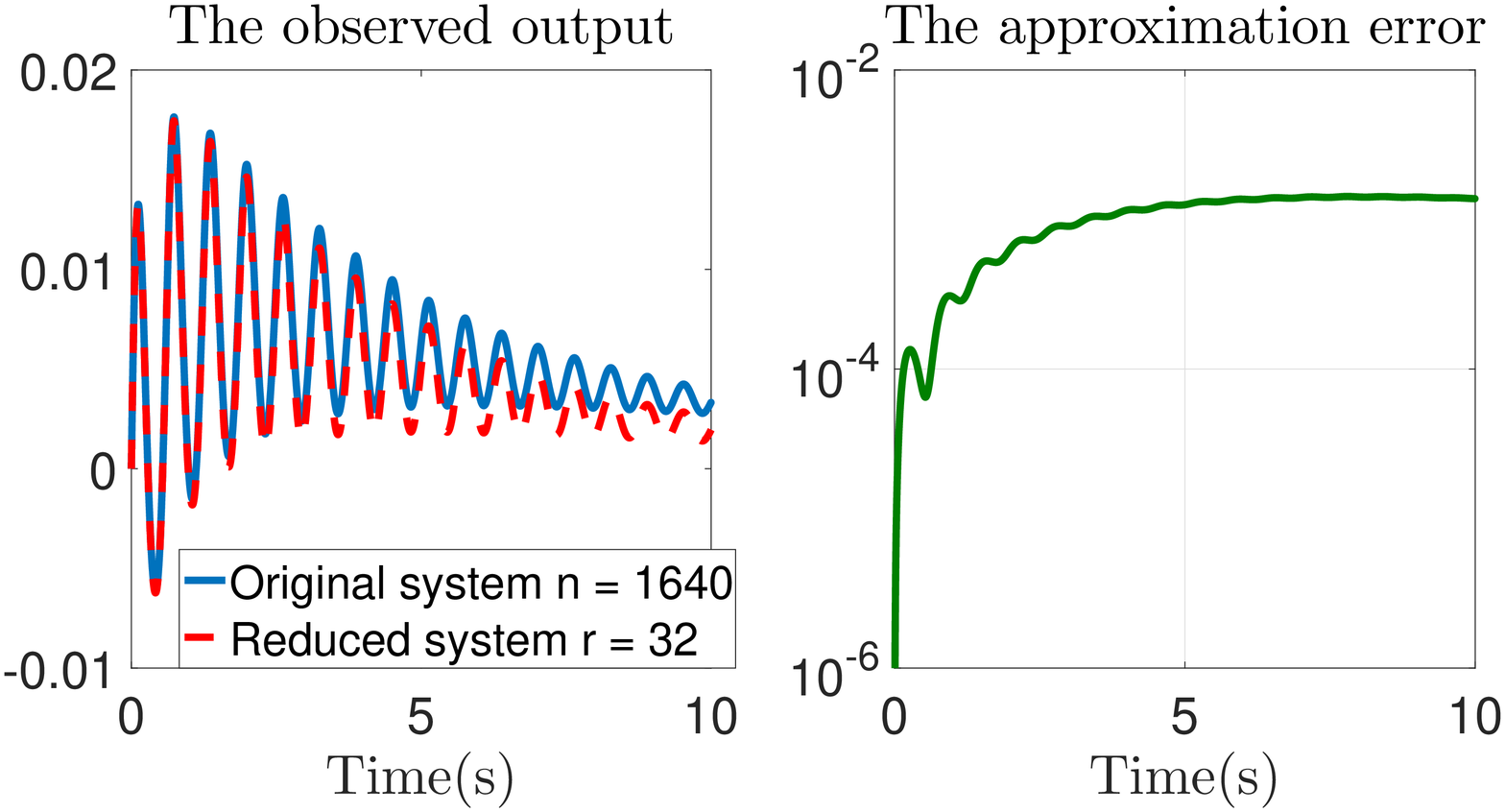} 
  		
  		\caption{Left plot: the observed outputs; 
  			right plot: the corresponding approximation error.}
  		\label{fig:GP3}
  	\end{center}
  	\vspace{-3mm}
  \end{figure}

  For the next experiment, choose the  tolerance value to be $\tau_r = 10^{-5}$. This corresponds to a truncation value of $r=32$. After computing the required matrices, notice that the D term is again numerically 0, while the norm of the matrix $\hat{\bF}$ slightly decreases to the value $6.9597e-04$. Perform numerical simulations and depict the two outputs and the approximation error in  Figure\;\ref{fig:GP3}. Observe that the approximation quality significantly improved, but there is still room for improvement.
  
\begin{figure}[!htb]
  	\vspace*{-2ex}
  	\begin{center}
  		\includegraphics[width = 0.9\textwidth]{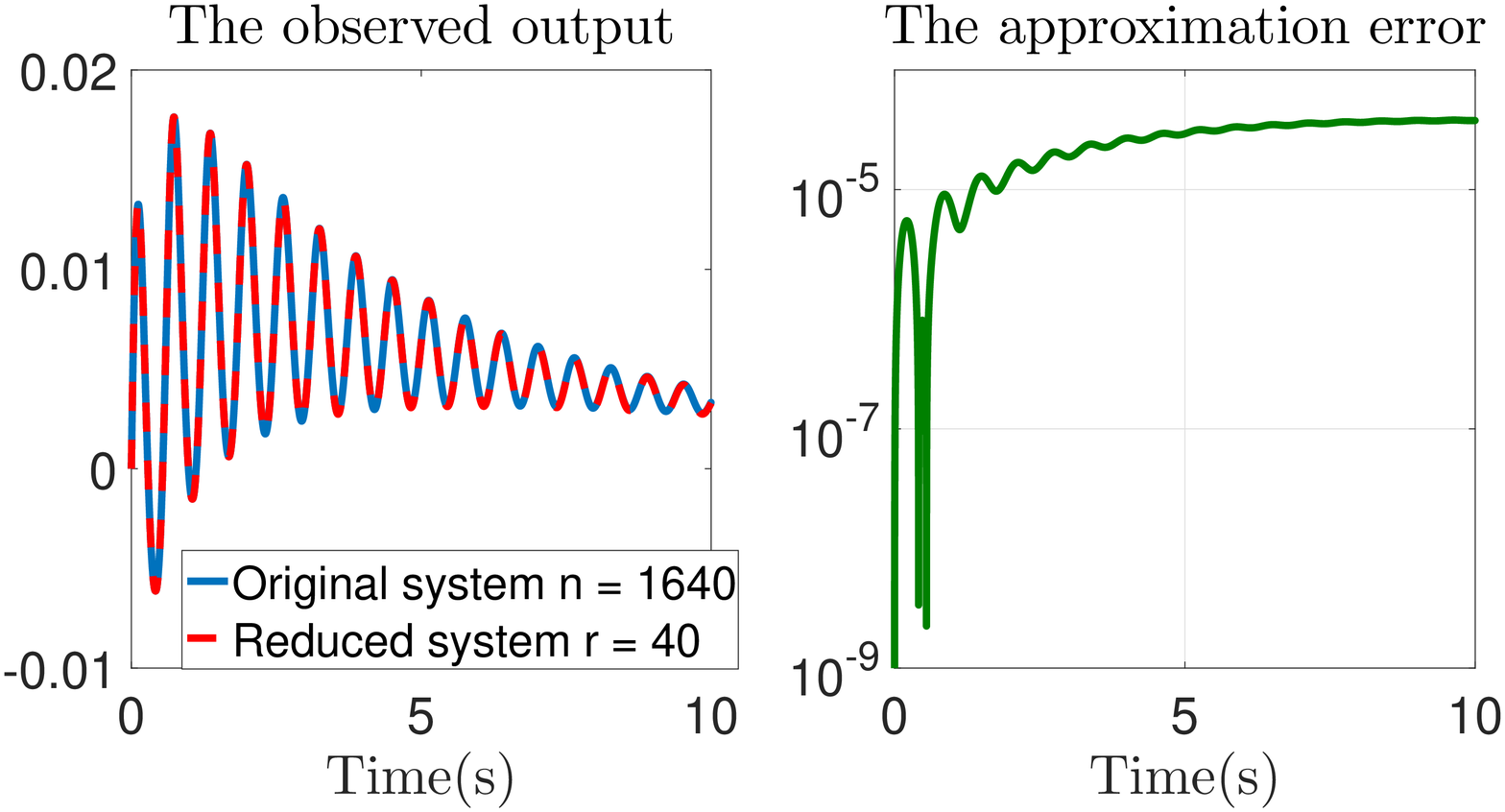} 
  		
  		\caption{Left plot: the observed outputs; 
  			right plot: the corresponding approximation error.}
  		\label{fig:GP4}
  	\end{center}
  	\vspace{-3mm}
  \end{figure}

  Finally, the tolerance value is chosen as $\tau_r = 10^{-6}$. For this particular choice, it follows that the truncation value is $r=40$. In this case, the output of the reduced-order model faithfully reproduces the original output, as it can be observed in Figure\;\ref{fig:GP4}. Note also that the approximation error stabilizes within the range $(10^{-4},10^{-5})$.
  
  \subsubsection{Experiment 2 - testing the trained models}\label{sec:Burgers_Exp2}

  In the second setup we perform time-domain simulations of the reduced-order models for different conditions than those used in the training stage, i.e. the time horizon is extended to $[0,15]s$ and two other control inputs are used. Moreover, we keep the truncation value to be $r=40$ (corresponding to tolerance $\tau_r = 10^{-6}$).
  
\begin{figure}[!htb]
  	\vspace*{-2ex}
  	\begin{center}
  		\includegraphics[width = 0.9\textwidth]{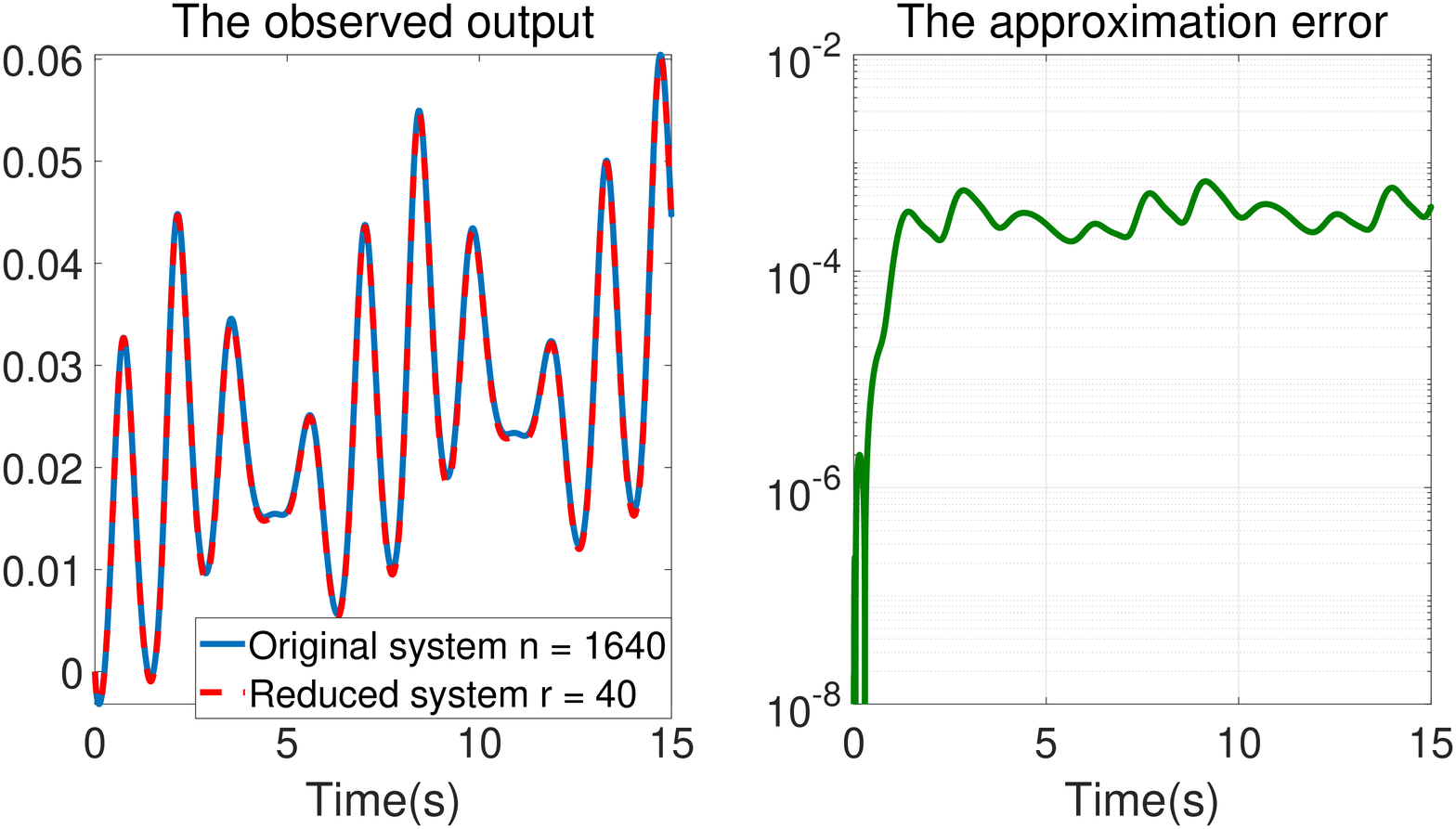} 
  		
  		\caption{Left plot: the observed outputs; 
  			right plot: the corresponding approximation error.}
  		\label{fig:GP5}
  	\end{center}
  	\vspace{-3mm}
  \end{figure}

  First, choose the testing control input to be $u_1(t) = \sin(4t)/4 - \cos(5t)/5$. The time domain simulations showing the observed outputs are depicted in Fig.\;\ref{fig:GP5}. Moreover, on the same figure, the magnitude of the approximation is presented. We observe that the output of the learned reduced model accurately approximates the output of the original system.
  
\begin{figure}[!htb]
  	\vspace*{-2ex}
  	\begin{center}
  		\includegraphics[width = 0.9\textwidth]{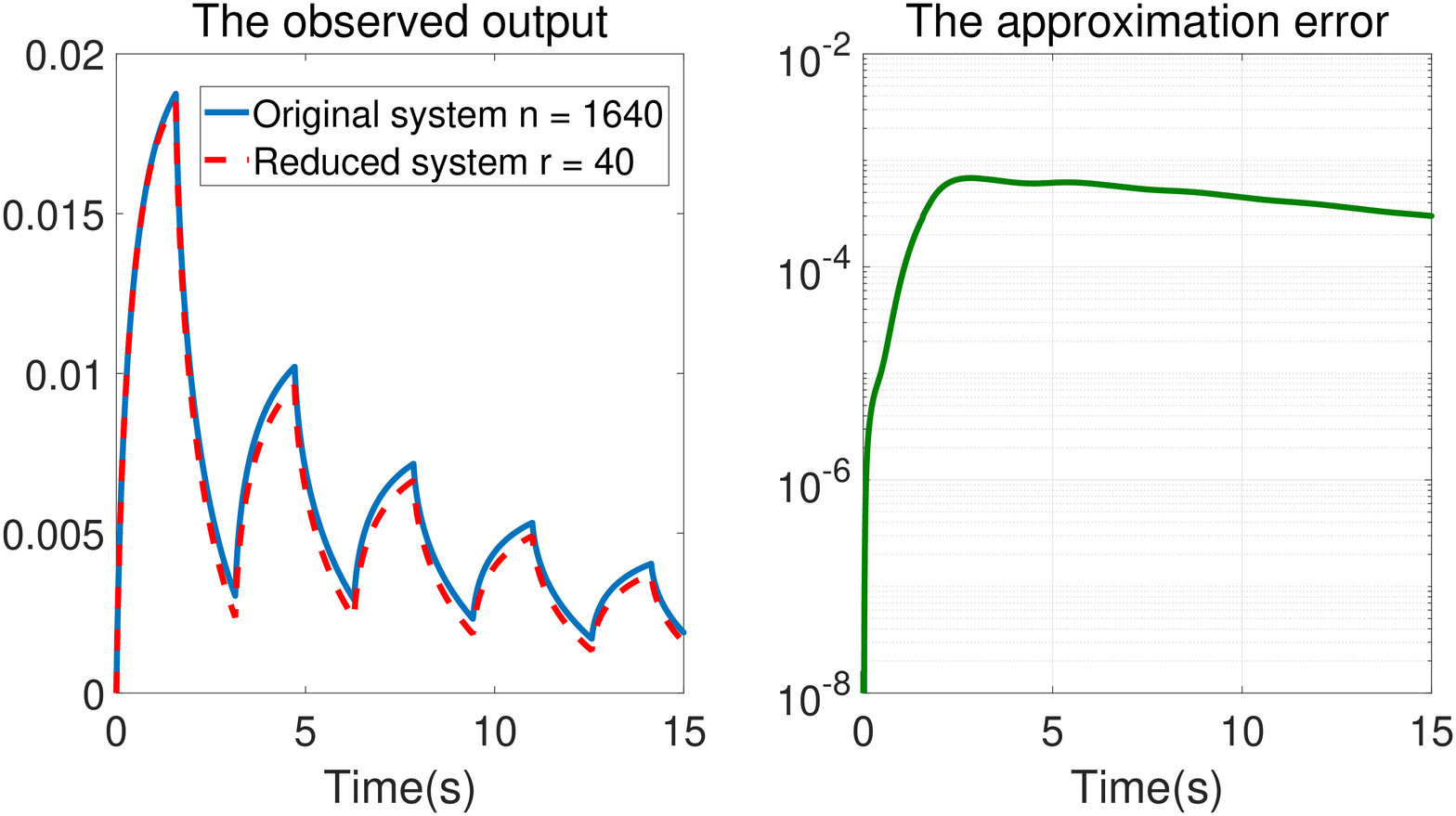} 
  		
  		\caption{Left plot: the observed outputs; 
  			right plot: the corresponding approximation error.}
  		\label{fig:GP6}
  	\end{center}
  	\vspace{-3mm}
  \end{figure}

  Afterwards, choose the testing control input to be $u_2(t) =$ \large $\frac{{\text square}(2t)}{5(t+1)}$ \normalsize. Note that ${\text square}(2t)$ is a square wave with period $\pi$. The time domain simulations showing the observed outputs are depicted in Fig.\;\ref{fig:GP6}. Moreover, on the same figure, the magnitude of the approximation is presented. We observe that the output of the learned reduced model does not approximate the output of the original system as well as in the previous experiments.
  
  \subsection{Coupled van der Pol oscillators}
  
  Consider the coupled van der Pol oscillators along a limit cycle example given in \cite{KS17}. The dynamics are characterized by the following six differential equations with linear and nonlinear (cubic) terms:
  \begin{align} \label{eq_vdp}
  \begin{split}
  \dot{x}_1 &= x_2, \\
  \dot{x}_2 &= -x_1 - \mu(x_1^2-1)x_2+a(x_3-x_1)+b(x_4-x_2), \\
  \dot{x}_3 &= x_4, \\
  \dot{x}_4 &= -x_3 - \mu(x_3^2-1)x_4+a(x_1-x_3)+b(x_2-x_4), \\
  &+ a(x_5-x_3) + b(x_6-x_4) +u, \\
  \dot{x}_5 &= x_6, \\
  \dot{x}_6 &= -x_5 - \mu(x_5^2-1)x_6+a(x_3-x_5)+b(x_4-x_6).
  \end{split}
  \end{align}
  Choose the output to be $y = x_3$. Hence the state-output equation is written as $ y = \bC \bx$ with $\bC = \left[ \begin{matrix}
  0 & 0 & 1 & 0 & 0 & 0 
  \end{matrix} \right]$ and $\bx = \left[ \begin{matrix}
  x_1 & x_2 & x_3 & x_4 & x_5 & x_6 
  \end{matrix} \right]^T$. Choose the parameters in (\ref{eq_vdp}) as follows: $\mu = 0.5, \ a = 0.5$ and  $b = 0.2$.

  Note that by introducing three additional surrogate states, e.g. $x_7 = x_1^2, x_8 = x_2^2$ and $x_9 = x_3^2$, one can rewrite the cubic nonlinear system in (\ref{eq_vdp}) of order $n=6$ as an order $n_q = 9$ quadratic-bilinear system. 
  
  Perform time-domain simulations of the cubic system of order $n=6$ and collect data from 500 snapshots using the explicit Euler method with step size $\Delta_t = 0.01$. The chosen time horizon is hence [0,5]s. The control input is a square-wave with period $\pi/5$ and amplitude 30, i.e., $u(t) = 30 \ \text square(10t)$.
  
  Compute the pseudo-inverse of matrix $\bOmega \in \mathbb{R}^{49 \times 500}$ and select as truncation value $p = 19$ (the 20th normalized singular value drops below machine precision).
  
  We compute a reduced-order quadratic-bilinear model of order $r=5$. We made this choice since the fifth normalized singular value of matrix $\bGamma \in \mathbb{R}^{7 \times 500}$ is  5.8651e-04  while the sixth is numerically 0, i.e., 3.5574e-16.
  We hence fit an order $r = 5$ quadratic-bilinear system that approximates the original order $n = 6$ cubic polynomial system.
  Note that the only non-zero feed-through quantity in the recovered state-output equation si given by $\hat{\bC} = \left[ \begin{matrix}
  -0.1067  &  0.5580  &  0.0797  & -0.4145 &  -0.7065
  \end{matrix} \right]$. 
  
\begin{figure}[!htb]
  	\vspace*{-2ex}
  	\begin{center}
  		\includegraphics[width = 0.9\textwidth]{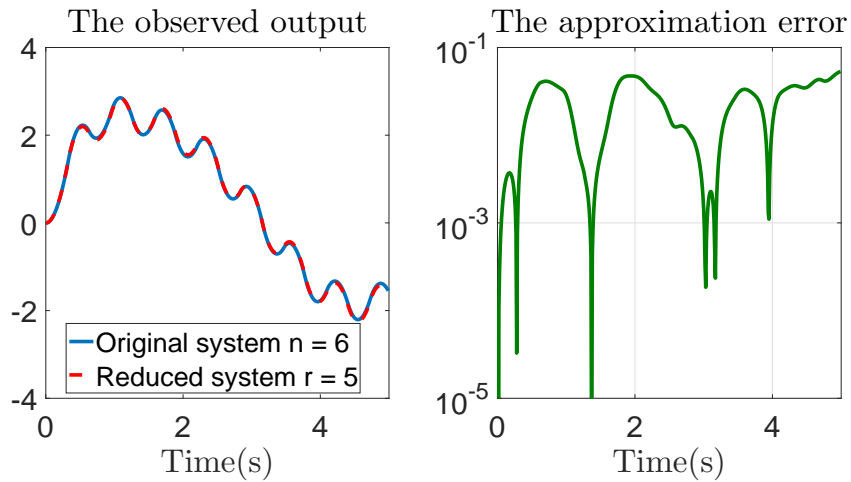} 
  		
  		\caption{Left plot: the observed outputs; 
  			Right plot: the corresponding approximation error.}
  		\label{fig:GP7}
  	\end{center}
  	\vspace{-3mm}
  \end{figure}

  Next, we perform time-domain simulations in the same manner as in Section \ref{sec:Burgers_Exp1}, i.e., by validating the reduced models on the training data. The results are depicted in Fig.\;\ref{fig:GP5}. One can observe that the two outputs match well. In this particular setup, it follows that the response of the 6th order cubic system (that can be equivalently written as a 9th order QB system) can be accurately approximated with the response of a 5th order QB system. The approximation error is presented in Fig.\;\ref{fig:GP7}.
  
\section{Conclusion}

In this paper, we have proposed extensions of the DMDc and ioDMD recently proposed methods. The philosophy is similar to that of the original methods, but instead of fitting discrete-time linear systems, we impose a more complex structure to the fitted models. More precisely, we fit bilinear or quadratic terms to augment the existing linear quantities (both in the differential and in the output equations). The numerical results presented were promising, and they have shown the strength of the method. Indeed, there is a clear trade-off to be made between approximation quality and the dimension of the fitted model. 

Nevertheless, this represents a first step towards extending DMD-type methods, and a more involved analysis of the method's advantages and disadvantages could represent an appealing future endeavor. Moreover, another contribution could be made by comparing the proposed methods in this work with the recently-introduced operator inference-type methods. For the quadratic-bilinear case, additional challenges arise when storing the large-scale matrices involved and also when computing the classical SVD for such big non-sparse matrices.


\section{Appendix}

\subsection{Computation of the reduced-order matrices for the quadratic-bilinear case}
\label{Appendix_QB}

In this section we present practical details for retrieving the system matrices in the case of the proposed procedure in Section \ref{sec:QB_DMD}. We solve the equation $\bGamma = \bG \bOmega$ for which the matrices are given as in (\ref{notQB}), i.e. the case without output observations. We again utilize an SVD, now performed on the matrix $\bOmega$, i.e.,
\begin{equation}
\bOmega = \bV \Si \bW^T \approx  \tilde{\bV} \tilde{\Si} \tilde{\bW}^T
\end{equation}
where the full scale and reduced SVD matrices have the following dimensions 
$$
\begin{cases} \bV \in \mathbb{R}^{(n^2+2n+1) \times (n^2+2n+1)}, \ \ \Si \in \mathbb{R}^{(n^2+2n+1) \times (m-1)}, \ \ \bW \in \mathbb{R}^{(m-1) \times (m-1)}, \\ \tilde{\bV} \in \mathbb{R}^{(n^2+2n+1) \times p}, \ \ \tilde{\Si} \in \mathbb{R}^{p \times p}, \ \ \tilde{\bW} \in \mathbb{R}^{(m-1) \times p}. \end{cases}
$$
The truncation index is denoted with $r$, and write as before $\bOmega^{\dagger} \approx \tilde{\bW} \tilde{\Si}^{-1} \tilde{\bV}^T.$

By splitting up the matrix $\bV^T$ as $\tilde{\bV}^T = [\tilde{\bV}_1^T \ \  \tilde{\bV}_2^T \ \ \tilde{\bV}_3^T \  \ \tilde{\bV}_4^T]$, with
$$
\tilde{\bV}_1,\tilde{\bV}_3  \in \mathbb{R}^{n \times r}, \ \ \tilde{\bV}_2 \in \mathbb{R}^{1 \times r}, \ \ \tilde{\bV}_4 \in \mathbb{R}^{n^2 \times r},
$$
recover the matrices
\begin{equation}
\overline{\bA} = \bX_s \tilde{\bW} \tilde{\Si}^{-1} \tilde{\bV}_1^T, \ \ \overline{\bB} = \bX_s \tilde{\bW} \tilde{\Si}^{-1} \tilde{\bV}_2^T, \ \ \overline{\bN} = \bX_s \tilde{\bW} \tilde{\Si}^{-1} \tilde{\bV}_3^T, \ \
\overline{\bQ} = \bX_s \tilde{\bW} \tilde{\Si}^{-1} \tilde{\bV}_4^T.
\end{equation}
Again, perform an additional SVD, e.g. $\bX_s  \approx  \hat{\bV} \tilde{\Si} \hat{\bW}^T$, where $\hat{\bV} \in \mathbb{R}^{(n+1) \times r}, \ \ \hat{\Si} \in \mathbb{R}^{r \times r}, \ \ \hat{\bV} \in \mathbb{R}^{(m-1) \times r}$. Using the transformation $\bx = \hat{\bV} \tilde{\bx}$, the following reduced-order approximations are computed:
\begin{align*}
\tilde{\bA} &= \hat{\bV}^T  \overline{\bA} \hat{\bV} = \hat{\bV}^T   \bX_s \tilde{\bW} \tilde{\Si}^{-1} \tilde{\bV}_1^T \hat{\bV} \in \mathbb{R}^{r \times r} , \\ \tilde{\bB}  &= \hat{\bV}^T  \overline{\bB}  = \hat{\bV}^T \bX_s \tilde{\bW} \tilde{\Si}^{-1} \tilde{\bV}_2^T \in \mathbb{R}^{r}, \\ \tilde{\bN} &= \hat{\bV}^T  \overline{\bA} \hat{\bV} = \hat{\bV}^T   \bX_s \tilde{\bW} \tilde{\Si}^{-1} \tilde{\bV}_3^T \hat{\bV} \in \mathbb{R}^{r \times r}, \\ \tilde{\bQ} &= \hat{\bV}^T  \overline{\bQ} ( \hat{\bV} \otimes \hat{\bV} ) = \hat{\bV}^T   \bX_s \tilde{\bW} \tilde{\Si}^{-1} \tilde{\bV}_2^T  (\hat{\bV}  \otimes \hat{\bV} ) \in \mathbb{R}^{r \times r^2}.
\end{align*}

 \bibliographystyle{spmpsci}
 \bibliography{ref_modred}

\end{document}